\documentclass[a4paper,12pt]{article}
\usepackage[utf8]{inputenc}
\usepackage{amsmath}
\usepackage{amssymb}
\usepackage{amsthm}
\usepackage{xcolor}

\theoremstyle{plain}
\newtheorem*{thm}{Theorem}%[section]

\newcommand{\Ric}{\mathrm{Ric}}

\def\M{{\sf M}}

\def\d{{\sf d}}
\def\m{{\sf m}}

\usepackage[notref,notcite,color]{showkeys}

\title{Corrigendum to\\ ``On the geometry of metric measure spaces. I.''\\ Acta Math.~196 (2006)}
\author{Karl-Theodor Sturm}

\date{\today}

\begin{document}

\maketitle
This is a corrigendum 
to Acta Math.~196 (2006) as well as
 to the follow-up publications
 JFA  259 (2010) and to JFA 260 (2011).
% where the same erroneous argument was used.
 
\begin{thm}
In the formulations of
the {\bf Tensorization Property} for
\begin{enumerate}
\item   the {\sf CD}$(K,\infty)$ condition, Proposition 4.16 in {\bf [1]} = Acta Math.~196 (2006) by St.,
\item the {\sf CD}$^*(K,N)$  condition, Theorem 4.1 in {\bf [2]} = JFA  259 (2010) by Kathrin Bacher  \& St., and
\item the  {\sf CD}$(K,N)$ condition, Theorem 1.1 in  {\bf [3]} = JFA 260 (2011) by Qintao Deng \& St.,
\end{enumerate}
the assumption \emph{non-branching} has to be replaced by the stronger assumption that all the metric measure spaces $(\M_i,\d_i,\m_i)$ for $i=1,\ldots,\ell$ %under consideration 
are 
\begin{itemize}
\item[\rm (a)] \emph{%normalized, 
infinitesimally Hilbertian spaces} (in the sense of {\bf [4,5]}) or
\item[\rm (b)] \emph{smooth Finslerian spaces} (in the sense of {\bf [6,7]}).
\end{itemize}
\end{thm}

Recall that a metric measure space $(\M,\d,\m)$ is called 
%\emph{normalized} if $\m(\M)=1$, and it is called
\emph{infinitesimally Hilbertian} if the Cheeger energy $\mathcal E$ on $\M$ is a quadratic functional on $L^2(\M,\m)$ or, in other words, if the associated heat flow is linear --- 
and note that \emph{Finslerian spaces} constitute the main class of metric measure spaces that are not infinitesimally Hilbertian.

\paragraph{The previous erroneous proofs} relied on an argument which can not be made rigorous
({\bf [1]}, p.~115, first two lines):
\emph{``According to Lemma 2.11 (iii), since M is non-branching and since the $\nu_{0,j}$ for $j=1, ..., n$
are mutually singular, also the $\eta_j$ for $j=1, ..., n$ must be mutually singular.''} %Hence, $\eta=\frac1n\sum_{j=1}^n \eta_j$ ... 

As thankfully   pointed out to the author first by G.~Savar\'e and then also by N.~Gigli,  the composition of the partial optimal transports is not necessarily an optimal transport of the composed marginals, and thus mutual singularity %disjointness of the  supports 
of the partial marginals (and  non-branching property) does not imply  mutual singularity %disjointness of the supports 
of the partial midpoints. % of ... are not necessarily disjoint. % and thus the (Boltzmann and Reny) entropies are not necessarily sub-additive. 

\bigskip
 
\paragraph{The corrected statements are derived as follows:} 
\subparagraph{(a) Infinitesimally Hilbertian spaces.}
\begin{itemize}
\item If all the $(\M_i,\d_i,\m_i)$  are \emph{infinitesimally Hilbertian spaces and satisfy the {\sf CD}$(K,\infty)$ condition} then so 
%is and 
does their tensor product $(\M,\d,\m)$, {\bf [8]}, Thms.~5.2 and 4.17. This is a consequence of the tensorization property for the Bakry-\'Emery condition  {\sf BE}$(K,\infty)$ and the equivalence between Eulerian and Lagrangian curvature-dimension conditions
%, i.e.~between 
{\sf BE}$(K,\infty)$  and {\sf CD}$(K,\infty)$.

(Note that 
the proof of the tensorization property  in {\bf [4]}, Theorem 6.13, is based on the incorrect Proposition 4.16 in {\bf [1]}. Also notice that in {\bf [8]} at some stage a tensorization argument from {\bf [4]} is used, but as noted there the lemma used was not involving curvature properties.)
%
%(Indeed, the original argument in {\bf [5]} only applied to non-branching normalized mm-spaces but it has been extended to non-branching $\sigma$-finite mm-spaces in {\bf [8]}, Thm.~7.6. The non-branching assumption -- in its weaker form of ``essentially non-branching'' -- was proven to follow from the other assumptions in {\bf [9]}, Cor.~1.2.)

\item 
If all the $(\M_i,\d_i,\m_i)$  are \emph{infinitesimally Hilbertian spaces and satisfy the {\sf CD}$^*(K,N)$ condition} then so %is and 
does their tensor product  $(\M,\d,\m)$, {\bf [9]}, Thm.~3.23.

\item 
If all the $(\M_i,\d_i,\m_i)$ are \emph{infinitesimally Hilbertian spaces and satisfy the {\sf CD}$(K,N)$ condition} then, in particular, they all satisfy the {\sf CD}$^*(K,N)$ condition. Thus by the previous assertion, their tensor product  $(\M,\d,\m)$,  is infinitesimally Hilbertian and satisfies the {\sf CD}$^*(K,N)$ condition. According to the Globalization Theorem of  {\bf [10]}, the latter implies that $(\M,\d,\m)$ also satisfies  the {\sf CD}$(K,N)$ condition.
(Indeed, the original argument in {\bf [10]} only applied to normalized mm-spaces but it is extended to $\sigma$-finite mm-spaces in {\bf [11]}.)
\end{itemize}
\subparagraph{(b) Smooth Finslerian spaces.}
 If all the $(\M_i,\d_i,\m_i)$ for $i=1,\ldots,\ell$ are \emph{smooth Finslerian spaces} (in the sense of {\bf [6,7]}) then so is their tensor product $(\M,\d,\m)=\bigotimes_{i=1}^\ell(\M_i,\d_i,\m_i)$.
For smooth Finslerian spaces, the {\sf CD}$(K,\infty)$ condition is equivalent to the \emph{weighted flag Ricci curvature} being bounded from below by $K$, {\bf [6]}, Thm.~1.2. By construction, the weighted flag Ricci curvature bound has the tensorization property.

Moreover, the {\sf CD}$^*(K,N)$ condition as well as the {\sf CD}$(K,N)$ condition are equivalent to the \emph{weighted flag $N$-Ricci curvature} being bounded from below by $K$, {\bf [6]}, Thm.~1.2.
 The tensorization property of  the weighted flag $N$-Ricci curvature bound follows as in the Riemannian case: if each of the spaces $\M_i$ with dimension $n_i$ and weight $V_i$ for $i=1,\ldots,\ell$ satisfies
 $$\Ric_{N_i}(\xi_i):=\Ric(\xi_i)+\mathrm{Hess } V_i(\xi,\xi)-\frac1{N_i-n_i}\langle\nabla V_i,\xi_i\rangle^2\ge K |\xi_i|^2\qquad\quad(\forall \xi_i\in T\M_i)$$ 
 then  the space $\M=\otimes_i \M_i$ with dimension $n=\sum_in_i$ and weight $V=\oplus_i V_i$ satisfies
\begin{align*} &\Ric(\xi)-\mathrm{Hess } V(\xi,\xi)-K |\xi|^2
=
\sum_i \Ric(\xi_i) -\sum_i \mathrm{Hess } V(\xi_i,\xi_i)-\sum_iK |\xi_i|^2\\
&\ge \sum_i \frac1{N_i-n_i}\, \langle\nabla V_i,\xi_i\rangle^2 \ge  \frac1{N-n}\, \langle\nabla V,\xi\rangle^2\qquad\quad (\forall \xi=(\xi_1,\ldots,\xi_\ell)\in T\M)
\end{align*}
 with 
$N=\sum_i N_i$. %  and $|\xi|^2=\sum_i|\xi_i|^2$.

\bigskip 

Indeed, in all these cases, the tensorization property in the Lagrangian picture is known only through its proof within the Eulerian picture.

%\subsubsection*{ad 1. Tensorization Property for the {\sf CD}$(K,\infty)$ condition}

%\begin{Prop} Let $(M_i,\d_i,m_i)$ for $i=1,\ldots.\ell$ be metric measure spaces and
%$$(M,\d,m)=\bigotimes_{i=1}^\ell(M_i,\d_i,m_i).$$
%Assume that $M$ is non-branching and compact. Then
%$$\underline{\text{\rm Curv}}(M,\d,m)=\inf_{i\in\{1,\ldots,\ell\}}\underline{\text{\rm Curv}}(M_i,\d_i,m_i).$$
%In other words, for any $K\in\R$,
%$$\forall i\in\{1,\ldots,\ell\}: \ (M_i,\d_i,m_i) \text{ satisfies }{\sf CD}(K,\infty) \quad\Longleftrightarrow\quad (M,\d,m) \text{ satisfies }{\sf CD}(K,\infty).$$
%\end{Prop}
%
\paragraph{References}
\begin{itemize}
\item[{[1]}] Karl-Theodor Sturm. ``On the geometry of metric measure spaces. I." Acta mathematica 196.1 (2006): 65-131.
%\item[{[2]}] Sturm, Karl-Theodor. ``On the geometry of metric measure spaces. II." Acta mathematica 196.1 (2006): 65-131.
 \item[{[2]}] Kathrin Bacher and Karl-Theodor Sturm. ``Localization and tensorization properties of the curvature-dimension condition for metric measure spaces." Journal of Functional Analysis 259.1 (2010): 28-56.
 
 \item[{[3]}]Qintao Deng and Karl-Theodor Sturm. ``Localization and tensorization properties of the curvature-dimension condition for metric measure spaces, II." Journal of Functional Analysis 260.12 (2011): 3718-3725.

    \item[{[4]}]  Luigi Ambrosio, Nicola Gigli, and Giuseppe Savar\'e. ``Metric measure spaces with Riemannian Ricci curvature bounded from below." Duke Mathematical Journal 163.7 (2014): 1405-1490.

     \item[{[5]}]  Nicola Gigli. ``On the differential structure of metric measure spaces and applications.''
  Memoirs of the AMS,  Volume 236, Number 1113 
(2015).
  
  \item[{[6]}] Shin-ichi Ohta. ``Finsler interpolation inequalities." Calculus of Variations and Partial Differential Equations 36.2 (2009): 211-249.
    
      \item[{[7]}] Shin-ichi Ohta and Karl-Theodor Sturm. ``Heat flow on Finsler manifolds." Communications on Pure and Applied Mathematics 62.10 (2009): 1386-1433.

             \item[{[8]}]  Luigi Ambrosio, Nicola Gigli, and Giuseppe Savar\'e. ``Bakry-\'Emery curvature-dimension condition and Riemannian Ricci curvature bounds.''
Annals of Probability 43.1 (2015): 339-404.
      
%       \item[{[8]}] Luigi Ambrosio,  Nicola Gigli, Andrea Mondino, and Tapio Rajala. ``Riemannian Ricci curvature lower bounds in metric measure spaces with $\sigma$-finite measure." Transactions of the American Mathematical Society 367.7 (2015): 4661-4701.
%       
%          \item[{[9]}] Tapio Rajala and Karl-Theodor Sturm. ``Non-branching geodesics and optimal maps in strong  CD(K,$\infty$)-spaces." Calculus of Variations and Partial Differential Equations 50.3 (2014): 831-846.
%          
  \item[{[9]}]  Matthias Erbar, Kazumasa Kuwada, and Karl-Theodor Sturm. ``On the equivalence of the entropic curvature-dimension condition and Bochner's inequality on metric measure spaces." Inventiones mathematicae 201.3 (2015): 993-1071.

\item[{[10]}] Fabio Cavalletti and Emanuel Milman. ``The globalization theorem for the curvature dimension condition." 
Inventiones mathematicae  226 (2021): 1-137.

 \item[{[11]}] Zhenhao Li:  ``On the Cavalletti-Milman globalization theorem." Annali di Matematica (2023). https://doi.org/10.1007/s10231-023-01352-9
\end{itemize}

\end{document}